\documentclass[11pt]{article}
\usepackage[margin=1in]{geometry}
\usepackage{amsmath,amssymb,amsthm,mathtools}
\usepackage{enumitem}
\usepackage{microtype}
\usepackage{aliascnt}
\usepackage{url}
\usepackage[hidelinks]{hyperref}
\usepackage[capitalise,noabbrev]{cleveref}
\usepackage[T1]{fontenc}
\usepackage{lmodern}

\numberwithin{equation}{section}
\newtheorem{theorem}{Theorem}[section]
\newaliascnt{proposition}{theorem}

\aliascntresetthe{proposition}
\newaliascnt{lemma}{theorem}
\newtheorem{lemma}[lemma]{Lemma}
\aliascntresetthe{lemma}
\newaliascnt{corollary}{theorem}

\aliascntresetthe{corollary}
\theoremstyle{definition}
\newaliascnt{definition}{theorem}
\newtheorem{definition}[definition]{Definition}
\aliascntresetthe{definition}
\theoremstyle{definition}
\newaliascnt{remark}{theorem}
\newtheorem{remark}[remark]{Remark}
\aliascntresetthe{remark}

\newcommand{\R}{\mathbb R}

\newcommand{\capm}{\operatorname{cap}_m}
\newcommand{\Rieszcap}{\operatorname{R}_{2m}}
\newcommand{\Capge}{\operatorname{Cap}^{\ge}_{m,\R^n}}

\newcommand{\supp}{\operatorname{supp}}

\title{Counterexamples for polyharmonic Dirichlet problems and $m$-harmonic capacity}
\author{Honggang Li}
\date{}

\usepackage{indentfirst}
\begin{document}
\maketitle
\begin{abstract}
In this paper, we construct counterexamples for polyharmonic Dirichlet problems
and $m$-harmonic capacity in $\mathbb{R}^n$. When $m=2$ and $n\ge9$, or $m\ge3$
and $n\ge2m+3$, we find a compact set $K\subset\mathbb{R}^n\setminus\{0\}$ whose
capacitary potential satisfies $U_K(0)>2$. Based on this result, we obtain a bounded connected domain $\Omega$
and $0\le f\in C_0^\infty(\Omega)$ so that the variational solution
of $(-\Delta)^m u=f$ is unbounded both above and below
along a nontangential sequence approaching a boundary
point satisfying Maz'ya's Wiener condition. The construction and
related capacitary arguments yield counterexamples to a pointwise
estimate for the Dirichlet Green function with a constant independent of
the domain, and show that $m$-harmonic capacity fails to be subadditive and can be strictly larger than the corresponding Riesz capacity. In these dimensional ranges, our results answer Maz'ya's Problems~32--35 \cite{MazyaProblems} and Problem~1 in \cite{Mazya1999Green}, as well as questions posed by Felli and Romani \cite{FelliRomani}.

\end{abstract}

\tableofcontents

\section{Introduction}\label{sec:introduction}

This paper is concerned with polyharmonic capacitary potentials, their applications to polyharmonic Dirichlet problems, and properties of higher-order capacities.

Let $m\ge 1$ be an integer and let $\Omega\subset\mathbb R^n$ be a bounded open set. We consider the polyharmonic Dirichlet problem
\begin{equation}
(-\Delta)^m u=f\quad\text{in }\Omega,\qquad u\in H_0^m(\Omega).
\label{eq:s1-1}
\end{equation}
A function $u\in H_0^m(\Omega)$ is a variational solution of \eqref{eq:s1-1} if
\[
\mathcal E_{m,\Omega}(u,v)=\int_\Omega fv\,dx
\qquad\text{for all }v\in H_0^m(\Omega),
\]
where
\[
\mathcal E_{m,\Omega}(u,v)=\sum_{|\alpha|=m}\frac{m!}{\alpha!}\int_\Omega D^\alpha u\,D^\alpha v\,dx.
\]
We write $\mathcal E_m=\mathcal E_{m,\mathbb R^n}$. Throughout the paper, all function spaces are real and $B_r=B(0,r)$.

For $f\in C_0^\infty(\Omega)$, the boundedness of variational solutions of \eqref{eq:s1-1} follows directly from the Sobolev embedding theorem when $n<2m$. In the critical dimension $n=2m$, Maz'ya \cite{Mazya1999Green} obtained logarithmic estimates for the Dirichlet Green function, which also imply the boundedness of such solutions. For $n>2m$, Maz'ya's pointwise Green-function estimates imply that the variational solution of \eqref{eq:s1-1} is bounded when $m=2$ and $5\le n\le7$ \cite{Mazya1979Green}, or when $m\ge3$ and $n=2m+1,\,2m+2$ \cite{Mazya1999Green}.

Under additional assumptions on the geometry of the domain, boundedness is also known in higher dimensions. For each fixed smooth bounded domain, this follows from classical elliptic regularity \cite{AgmonDouglisNirenberg1959}. Mayboroda and Maz'ya \cite{MayborodaMazya2008Convex} established the boundedness of variational solutions to the biharmonic Dirichlet problem on arbitrary bounded convex domains in all dimensions; see \cite{MayborodaMazya2009,MayborodaMazya2014,MayborodaMazya2018} for further regularity results for polyharmonic equations. More recently, Lemenant and Mougenot \cite{LemenantMougenot2026} proved boundary regularity for polyharmonic Dirichlet problems on Reifenberg-flat domains with sufficiently small flatness parameters. For arbitrary bounded open sets, however, the question remained open in higher dimensions. Maz'ya stated it as Problem~33 in~\cite{MazyaProblems}: are variational solutions of \eqref{eq:s1-1} bounded when $m=2$ and $n\ge8$, or when $m\ge3$ and $n\ge2m+3$?

A related question concerns the behavior of variational solutions near a
boundary point.  Assume $n>2m$. For a compact set $F\subset\mathbb{R}^n$, let
\[
\mathcal C_F
=\{v\in C_0^\infty(\mathbb R^n):
v=1\text{ in a neighborhood of }F\},
\qquad
\capm(F)=\inf_{v\in\mathcal C_F}\mathcal E_m(v,v).
\]
In 1924, Wiener \cite{Wiener1924} characterized regular boundary points for the
Laplacian in terms of harmonic capacity. For $(-\Delta)^m$, a point
$0\in\partial\Omega$ is called regular if the variational solution of
\eqref{eq:s1-1} satisfies $u(x)\to 0$ as $x\to 0$ in $\Omega$ for every
$f\in C_0^\infty(\Omega)$. When $m=2$ and $n=5,6,7$, or $m\ge 3$ and
$n=2m+1,2m+2$, Maz'ya \cite[Theorem~2]{MazyaWiener2002} proved that $0$ is regular
if and only if
\begin{equation}\label{eq:wiener}
\int_0^1 \frac{\operatorname{cap}_m(\overline{B_r}\setminus\Omega)}
{\operatorname{cap}_m(\overline{B_r})}\,\frac{dr}{r}=\infty.
\end{equation}
Here we use the equivalence of $\operatorname{cap}_m$
and Bessel capacity on compact subsets of $B_1$;
see \cite[Section~2]{MazyaWiener2002}. For $m=2$ and
$n\ge 8$, Maz'ya \cite[Section~10, Theorem~3]{MazyaWiener2002} also showed that
\eqref{eq:wiener} implies nontangential convergence to zero at the
vertex of certain inner cusps; see \cite{MazyaBoundary2018} for further results.
In \cite[Problem~1]{Mazya1999Green}, he asked whether \eqref{eq:wiener} is
sufficient for regularity when $m=2$ and $n\ge 8$, or $m\ge 3$ and
$n\ge 2m+3$.

We give negative answers to both questions in all of these cases except
$m=2$ and $n=8$. Our counterexamples are based on polyharmonic capacitary potentials.
Let $\mathring W^{m,2}(\mathbb R^n)$ be the completion of
$C_0^\infty(\mathbb R^n)$ under the norm
$\mathcal E_m(v,v)^{1/2}$.
The polyharmonic capacitary potential $U_F$ is the unique minimizer
satisfying
\[
U_F\in\overline{\mathcal C_F}^{\,\mathring W^{m,2}},
\qquad
\mathcal E_m(U_F,U_F)=\capm(F).
\]
For $m=1$, capacitary potentials take values in $[0,1]$,
and the corresponding capacity is subadditive; see
\cite{TaylorPDEII2011}.
For $m\ge2$, however, these properties do not follow from
the arguments used in the case $m=1$.

Although the bound $U_F\le1$ need not hold for higher-order
capacitary potentials, Maz'ya \cite[Section~4]{MazyaWiener2002}
proved that, when $m=2$ and $n=5,6,7$, or when $m\ge3$ and
$n=2m+1,2m+2$, every compact set $F$ of positive $m$-harmonic
capacity satisfies
\[
0<U_F(x)<2,\qquad x\in\mathbb R^n\setminus F.
\]

In Problem~34 of \cite{MazyaProblems}, Maz'ya asked whether the
same bound holds in higher dimensions and whether the constant $2$
is optimal. Problem~35 concerns the subadditivity of $m$-harmonic
capacity. We show that the proposed higher-dimensional bound and
subadditivity both fail in the dimensional ranges specified below.

Our main results concern the following two cases:
\begin{equation}\label{eq:s1-5}
m=2,\qquad n\ge9,
\end{equation}
or
\begin{equation}\label{eq:s1-6}
m\ge3,\qquad n\ge2m+3.
\end{equation}
We first state the two main results. The first disproves the higher-dimensional
bound in Problem~34. The second uses the same construction to answer
Problem~33 and the question concerning the sufficiency of \eqref{eq:wiener}.

\begin{theorem}\label{thm:main-overshoot}
Suppose that either \eqref{eq:s1-5} or \eqref{eq:s1-6} holds. Then there exists a
compact set $K\Subset\R^n\setminus\{0\}$, which is a finite union of
pairwise disjoint closed balls, such that
\[
       U_K(0)>2.
\]
\end{theorem}
We also show that the constant $2$ in Maz'ya's estimate
cannot be replaced by a smaller constant independent of
$m$ and $n$; see \Cref{rem:uniform-sharpness}.

Based on the capacitary construction in \Cref{thm:main-overshoot},
we construct a variational solution that is unbounded both above and below
near a boundary point satisfying the Wiener condition.
\begin{theorem}\label{thm:main-unbounded}
Suppose that either \eqref{eq:s1-5} or \eqref{eq:s1-6} holds.
Then there exist a bounded connected domain
$\Omega\subset\mathbb{R}^n$ with $0\in\partial\Omega$
satisfying the Wiener condition \eqref{eq:wiener}, and a
function $f\in C_0^\infty(\Omega)$ with $f\ge 0$ and
$f\not\equiv 0$, with the following property.
For the variational solution $u\in H_0^m(\Omega)$ of
$(-\Delta)^m u=f$, there exist a unit vector
$e\in\mathbb{S}^{n-1}$, a constant $c>0$, and a sequence
$t_j\downarrow 0$ such that $x_j=t_j e\in\Omega$ and
\[
\operatorname{dist}(x_j,\partial\Omega)\ge c|x_j|,
\qquad
u(x_{2j})\longrightarrow+\infty,
\qquad
u(x_{2j+1})\longrightarrow-\infty.
\]
In particular, $\operatorname*{ess\,sup}_\Omega u=+\infty$
and $\operatorname*{ess\,inf}_\Omega u=-\infty$.
\end{theorem}

\Cref{thm:main-unbounded} gives negative answers to Problem~33
in \cite{MazyaProblems} and Problem~1 in \cite{Mazya1999Green}
under \eqref{eq:s1-5} or \eqref{eq:s1-6}. This also shows that positivity of the
fundamental solution of $(-\Delta)^m$ alone is not sufficient for the
Wiener criterion in \cite[Theorem~2]{MazyaWiener2002}.

Our construction also has consequences for pointwise estimates
of Dirichlet Green functions. In Problem~2 of \cite{Mazya1999Green}
and Problem~32 of \cite{MazyaProblems}, Maz'ya asked whether the Green function for
\eqref{eq:s1-1} satisfies
\[
|G_{\Omega}(x,y)|\le C(m,n)|x-y|^{2m-n}
\]
for all distinct $x,y\in\Omega$, where $C(m,n)$ is independent of
$\Omega$. Such an estimate is known on arbitrary bounded open sets when
$m=2$ and $n=5,6,7$, or when $m\ge3$ and $n=2m+1,2m+2$
\cite{Mazya1979Green,Mazya1999Green}. On the other hand, for each fixed smooth
bounded domain, an estimate of this form
holds with a constant that may depend on
the domain \cite{Krasovskii1967,DallAcquaSweers2004}.
The following theorem shows that, under
\eqref{eq:s1-5} or \eqref{eq:s1-6}, no constant
depending only on $m$ and $n$ suffices,
even among smooth bounded connected domains.

\begin{theorem}\label{thm:main-green}
Suppose that either \eqref{eq:s1-5} or \eqref{eq:s1-6} holds.
Then there exist smooth bounded connected
domains $\Omega_j^\pm\subset B_4$ with
$0\in\Omega_j^\pm$ and points
$y_j^\pm\in B_3\setminus\overline{B_2}$
such that
\begin{equation}
G_{\Omega_j^+}(0,y_j^+)\longrightarrow+\infty,
\qquad
G_{\Omega_j^-}(0,y_j^-)\longrightarrow-\infty.
\label{eq:s1-10}
\end{equation}

Consequently, for every $M>0$, there exist
a smooth bounded connected domain
$\Omega\subset\R^n$ and points $x,y\in\Omega$ with $x\ne y$
such that
\begin{equation}
|x-y|^{n-2m}|G_{\Omega}(x,y)|>M.
\label{eq:s1-11}
\end{equation}
\end{theorem}

The construction in \Cref{thm:main-overshoot}, together with a capacitary
perturbation argument, also yields two further results concerning
higher-order capacities. The first gives a stronger negative answer
to Maz'ya's Problem~35 on the subadditivity of $m$-harmonic capacity.
\begin{theorem}\label{thm:main-nonsubadd}
Suppose that either \eqref{eq:s1-5} or \eqref{eq:s1-6} holds.
For every $M>1$, there exist disjoint
compact sets $E,F\subset\R^n$, with $E$ a finite union of pairwise
disjoint closed balls and $F$ a single closed ball, such that
\[
       \capm(E\cup F)>\capm(E)+M\capm(F).
\]
\end{theorem}

The second concerns the relationship between two definitions of
higher-order capacity studied by Felli and Romani
\cite{FelliRomani}. For a compact set $F\subset\R^n$,
they introduced the capacity \cite[(2.8)]{FelliRomani}
\[
\Capge(F)
=\inf\left\{
\mathcal E_m(v,v):
v\in\mathring W^{m,2}(\R^n),\
v\ge1\text{ a.e. on }F
\right\}.
\]
Felli and Romani \cite[Section~5]{FelliRomani}
asked whether the two capacities coincide when $m\ge 2$,
or at least whether
\[
\capm(\varepsilon F)\sim\Capge(\varepsilon F)
\quad\text{as }\varepsilon\downarrow0
\]
for compact sets $F$ of positive capacity.

We give negative answers to both questions
by comparing these capacities with the
Riesz capacity \cite[(5.12)]{MazyaProblems},
defined by
\[
\Rieszcap(F)
=\inf\left\{
\mathcal E_m(v,v):
v\in C_0^\infty(\R^n),\
v\ge1\text{ on }F
\right\}.
\]
The Riesz capacity is subadditive
\cite[Section~5.14]{MazyaProblems}.

\begin{theorem}\label{thm:main-capacity-gap}
Suppose that either \eqref{eq:s1-5} or \eqref{eq:s1-6} holds. Then there exists a compact
set $H$, which is a finite union of pairwise disjoint closed balls,
such that
\begin{equation}
       \Capge(H)\le\Rieszcap(H)<\capm(H).
       \label{eq:s1-13}
\end{equation}
Moreover, for $H_\varepsilon=\varepsilon H$ and any $\varepsilon>0$,
\begin{equation}
 \frac{\Capge(H_\varepsilon)}{\capm(H_\varepsilon)}
 =\frac{\Capge(H)}{\capm(H)}<1.
 \label{eq:s1-14}
\end{equation}
\end{theorem}

We briefly describe the proofs of the main
results. To prove \Cref{thm:main-overshoot}, we first consider
functions whose trace on the unit sphere is 1.
Let $W$ be the function that minimizes
$\mathcal E_{m}(v,v)$ among them.
By the sharp trace inequality,
we obtain an explicit formula for $W$,
which shows that $W(0)>2$
under \eqref{eq:s1-5} or \eqref{eq:s1-6}.
We then construct a finite union $K_N$
of pairwise disjoint closed balls whose centers lie on the unit sphere and prove that
\[
U_{K_N}(0)\longrightarrow W(0)
\quad\text{as }N\to\infty.
\]
This gives the counterexample in \Cref{thm:main-overshoot}.

The proofs of \Cref{thm:main-unbounded,thm:main-green} use the set $K$ constructed in \Cref{thm:main-overshoot}.
By successively removing the sets $r_jK$ from a ball, where $r_j\downarrow0$, we construct a sequence of domains $\Omega_N$.
Let $u_N$ be the corresponding variational solution for a fixed nonnegative right-hand side $f$.
For each $N$, with $r_1,\ldots,r_{N-1}$ fixed and $u_{N-1}(0)\ne0$, we have
\[
\frac{u_N(0)}{u_{N-1}(0)}
\longrightarrow 1-U_K(0)<-1
\quad\text{as }r_N\downarrow0.
\]
By choosing the radii $r_j$ suitably, we obtain a variational solution
unbounded above and below along a nontangential sequence approaching the origin. Each removed set $r_jK$ contributes a fixed positive lower
bound to the Wiener integral over a corresponding interval of radii,
so that \eqref{eq:wiener} holds at the origin.
The Green representation formula for the solutions on $\Omega_N$ gives \Cref{thm:main-green}.

To prove \Cref{thm:main-nonsubadd,thm:main-capacity-gap}, we consider
a compact set $K$ and a small set
$L_r=x_0+rL$, where $x_0\notin K$.
We prove that
\[
\operatorname{cap}_m(K\cup L_r)
-\operatorname{cap}_m(K)
=
|1-U_K(x_0)|^2\operatorname{cap}_m(L_r)
+o(\operatorname{cap}_m(L_r))
\quad\text{as }r\downarrow0.
\]
Using this formula, we construct a finite union of pairwise disjoint closed balls
$E$ such that $|1-U_E(0)|$ is arbitrarily large.
For any $M>1$, choose $E$ with $|1-U_E(0)|^2>M$ and let
$L=\overline{B_1}$. For sufficiently small $r>0$, this gives
\[
\operatorname{cap}_m(E\cup rL)
>
\operatorname{cap}_m(E)+M\operatorname{cap}_m(rL),
\]
which yields \Cref{thm:main-nonsubadd}.
Using the subadditivity of Riesz capacity
and the scaling properties of the capacities,
we also obtain \Cref{thm:main-capacity-gap}.

The paper is organized as follows.
\Cref{sec:framework} contains the preliminaries.
\Cref{sec:sieve} proves \Cref{thm:main-overshoot}.
\Cref{sec:amplifier} contains the proofs of \Cref{thm:main-unbounded,thm:main-green}.
\Cref{thm:main-nonsubadd,thm:main-capacity-gap}
are proved in \Cref{sec:potential-amplifier}.

\section{Preliminaries}\label{sec:framework}
For an open set $G\subset\mathbb R^n$,
 $\mathring W^{m,2}(G)$ is the completion of $C_0^\infty(G)$ in
$\mathring W^{m,2}(\mathbb R^n)$,
with functions extended by zero.

\subsection{Mosco convergence}\label{sec:mosco}

We use Mosco convergence for closed subspaces of a Hilbert space $H$. Recall that $V_j$ converges to $V$ in this sense if
\begin{enumerate}
\item[$(\mathrm{M1})$] for every $v\in V$, there exist $v_j\in V_j$ such that $v_j\to v$ in $H$;
\item[$(\mathrm{M2})$] if $j_k\to\infty$, $v_k\in V_{j_k}$, and $v_k\rightharpoonup v$ in $H$, then $v\in V$.
\end{enumerate}
See \cite[Definition 1.1, pp.~519--520]{Mosco1969}.
For a family $\{V_r\}_{r>0}$ of closed subspaces
of $H$, we say that $V_r$ converges to $V$
in the sense of Mosco as $r\to0^+$
if the following two conditions are satisfied:
\begin{enumerate}
\item[$(\mathrm{M1}_r)$] For every $v\in V$, there exist $v_r\in V_r$ such that $v_r\to v$ in $H$ as $r\to0^+$;
\item[$(\mathrm{M2}_r)$] For every sequence $r_j\to 0^+$, if $v_j\in V_{r_j}$ and $v_j\rightharpoonup v$ in $H$, then $v\in V$.
\end{enumerate}
The Mosco convergence of $V_r$ to $V$
as $r\to0^+$ is equivalent to
the Mosco convergence of $(V_{r_j})$ to $V$
for every sequence $r_j\to0^+$.

\begin{lemma}\label{lem:mosco-minimizers}
Let $H$ be a Hilbert space. Let $V_j, V$ be closed subspaces of $H$ such that $V_j$ converges to $V$ in the sense of Mosco, and suppose that $b_j\to b$ in $H$. Let $z_j$ and $z$ minimize the norm on $b_j+V_j$ and $b+V$, respectively. Then
$$z_j\longrightarrow z\qquad\text{in }H.$$
\end{lemma}

\begin{proof}
Since $b_j\to b$ in $H$ and $V_j$ converges to $V$ in the sense of Mosco, the affine spaces $b_j+V_j$ converge to $b+V$ in the sense of Mosco. The conclusion follows from \cite[Corollary of Theorem~A, p.~542]{Mosco1969}.
\end{proof}

The conclusion of \Cref{lem:mosco-minimizers} also holds for families indexed by $r>0$ as $r\downarrow 0$, under the corresponding assumptions.

\subsection{The \texorpdfstring{$m$}{m}-harmonic capacity}\label{sec:m-harmonic-capacity}
Let $K\subset\mathbb R^n$ be a compact set. Recall that
\[
\mathcal C_K=\{v\in C_0^\infty(\mathbb R^n):
v=1\text{ in a neighborhood of }K\}.
\]
Choose $\eta_K\in\mathcal C_K$ and set
$$\mathcal A_K=\eta_K+\mathring W^{m,2}(\mathbb R^n\setminus K).$$
Since $\eta_1-\eta_2\in
C_0^\infty(\mathbb R^n\setminus K)
\subset\mathring W^{m,2}(\mathbb R^n\setminus K)$
for two such functions $\eta_1$ and $\eta_2$,
the space $\mathcal A_K$ is independent
of the choice of $\eta_K$.

\begin{definition}\label{def:capacity}
Define
\begin{equation}\label{eq:s2-11}
\operatorname{cap}_m(K)=\inf_{v\in\mathcal A_K}\mathcal E_m(v,v).
\end{equation}
\end{definition}

Since $\eta_K\in\mathcal C_K$,
we have
$\overline{\mathcal C_K}^{\,\mathring W^{m,2}}
=\mathcal A_K$.
Thus \Cref{def:capacity} coincides with Maz'ya's definition
of $m$-harmonic capacity \cite[Section~5.13, (5.11)]{MazyaProblems}.
By the Hilbert projection theorem, \eqref{eq:s2-11} has a
unique minimizer $U_K$. It satisfies
$$\mathcal E_m(U_K,\varphi)=0\qquad\text{for all }\varphi\in\mathring W^{m,2}(\mathbb R^n\setminus K).$$
In particular, $U_K$ is $m$-harmonic in $\mathbb R^n\setminus K$, and hence $U_K\in C^\infty(\mathbb R^n\setminus K)$. By rescaling and the uniqueness of the minimizer,
for every $r>0$ we have
$$
\operatorname{cap}_m(rK)
=r^{n-2m}\operatorname{cap}_m(K),
\qquad
U_{rK}(x)=U_K(x/r)
\quad (x\in\mathbb R^n\setminus rK).
$$

Moreover, the definitions of the three capacities yield
\begin{equation}\label{eq:capacity-comparison}
\operatorname{Cap}_{m,\mathbb R^n}^{\ge}(K)\le\operatorname{R}_{2m}(K)\le\operatorname{cap}_m(K).
\end{equation}

The following simple observation will be useful.
\begin{lemma}\label{lem:unit-neighborhood}
If $v\in\mathring W^{m,2}(\mathbb R^n)$ and $v=1$ almost everywhere in an open neighborhood of $K$, then $v\in\mathcal A_K$.
\end{lemma}

\begin{proof}

Let $O$ be an open neighborhood of $K$ such that
$v=1$ a.e.\ in $O$. Choose an open set $U$ with
$K\subset U\Subset O$ and a function $\eta_K\in C_0^\infty(O)$
such that $\eta_K=1$ in $U$, and set $w=v-\eta_K$.
Then $w\in\mathring W^{m,2}(\mathbb R^n)$
and $w=0$ a.e.\ in $U$.

Choose $\chi\in\mathcal C_K$ with $\supp\chi\subset U$, and let
$w_j\in C_0^\infty(\mathbb R^n)$ converge to $w$
in $\mathring W^{m,2}(\mathbb R^n)$.
By the homogeneous Sobolev inequalities
and standard cutoff estimates,
\[
(1-\chi)w_j\longrightarrow (1-\chi)w=w
\quad\text{in }\mathring W^{m,2}(\mathbb R^n).
\]
Since $(1-\chi)w_j\in
C_0^\infty(\mathbb R^n\setminus K)$,
we have
$w\in\mathring W^{m,2}(\mathbb R^n\setminus K)$.
Hence $v=\eta_K+w\in\mathcal A_K$.
\end{proof}

\subsection{A minimization problem on the sphere}\label{sec:membrane}

Let $S=\mathbb S^{n-1}$, and let $\sigma$ be the normalized surface measure on $S$. Consider
$$\mathcal M_S=\{v\in\mathring W^{m,2}(\mathbb R^n):\operatorname{Tr}_S v=1\}.$$
By the Hilbert projection theorem, there exists
a unique $W\in\mathcal M_S$ such that
$\mathcal E_m(W,W)=\inf\limits_{v\in\mathcal M_S}
\mathcal E_m(v,v)$.
It follows from \cite[Theorem~1.1]{BezMachiharaSugimoto2016} that
$$W(x)=\frac{P(x)}{I_{n,m}},\qquad P(x)=\int_S |x-\omega|^{2m-n}\,d\sigma(\omega),$$
where
$$I_{n,m}=\int_S |e-\omega|^{2m-n}\,d\sigma(\omega),\qquad e\in S.$$
For $x\ne\omega$ and $|\alpha|\le m$,
\begin{equation}\label{eq:s2-21}
|D_x^\alpha(|x-\omega|^{2m-n})|\le C|x-\omega|^{2m-n-|\alpha|}.
\end{equation}
Since $m\ge 2$ and $2m-|\alpha|-2\ge 0$, the right-hand side of \eqref{eq:s2-21} is uniformly integrable on $S$ for $x$ in a fixed neighborhood of $S$. Hence $D^\alpha W$ is bounded in a fixed neighborhood of $S$ for all $|\alpha|\le m$. In particular, $W$ is $C^1$ near $S$ and $W\in C^\infty(\mathbb R^n\setminus S)$.

By \cite[(3.12)]{BrauchartGrabner2015}, we have
$$I_{n,m}=\frac{2^{2m-2}\Gamma(n/2)\Gamma(m-\frac{1}{2})}{\sqrt{\pi}\,\Gamma(n/2+m-1)}.$$
Since $P(0)=1$,
\begin{equation}\label{eq:s2-23}
W(0)=\frac{\sqrt{\pi}\,\Gamma(n/2+m-1)}{2^{2m-2}\Gamma(n/2)\Gamma(m-\frac{1}{2})}.
\end{equation}
For $m\ge 2$ and $n>2m$,
\begin{equation}\label{eq:s2-24}
W(0)>2 \quad\text{if and only if}\quad
\begin{cases}
m=2, & n\ge 9,\\
m\ge 3, & n\ge 2m+3.
\end{cases}
\end{equation}
To see this, formula \eqref{eq:s2-23} gives
$W(0)=n/4$ when $m=2$.
For fixed $m\ge 3$, the same formula and
$\Gamma(z+1)=z\Gamma(z)$ give
$$
W(0)=\frac{\sqrt{\pi}}
{2^{2m-2}\Gamma(m-\frac{1}{2})}
\prod_{j=0}^{m-2}\left(\frac{n}{2}+j\right),
$$
so $W(0)$ is strictly increasing in $n$.
At $n=2m+2$, we have $W(0)=2-1/m<2$.
At $n=2m+3$, set $J_m=I_{2m+3,m}$. Then
$$
J_3=\frac{16}{33}<\frac{1}{2}
\qquad\text{and}\qquad
\frac{J_{m+1}}{J_m}
=\frac{4(2m-1)(2m+3)}
{(4m+1)(4m+3)}<1
\quad\text{for }m\ge 3.
$$
Thus $W(0)>2$ at $n=2m+3$,
and the monotonicity in $n$ yields
\eqref{eq:s2-24}.

\subsection{Partitions and estimates on the sphere}\label{sec:sphere-estimates}

Recall that $S=\mathbb S^{n-1}$, and set $d=n-1$.
Let $d_S$ denote the geodesic distance on $S$. For any $x,y\in S$,
\[
|x-y|\le d_S(x,y)\le
\frac{\pi}{2}|x-y|.
\]
For $\omega\in S$ and $r>0$, set
\[
B_S(\omega,r)
=\{x\in S:d_S(\omega,x)<r\}.
\]

\begin{lemma}\label{lem:spherical-partition}
There exist constants $c_0,C_0>0$ such that, for every sufficiently large $N$, there are $\sigma$-measurable sets $Q_{j,N}\subset S$ and points $\omega_{j,N}\in Q_{j,N}$, $j=1,\dots,N$, satisfying
$$
S=\bigcup_{j=1}^N Q_{j,N},\qquad
Q_{i,N}\cap Q_{j,N}=\varnothing\quad\text{for }i\ne j,\qquad
\sigma(Q_{j,N})=\frac{1}{N},
$$
and
$$
\operatorname{diam}_{d_S}(Q_{j,N})
=\sup_{x,y\in Q_{j,N}}d_S(x,y)
\le 2C_0N^{-1/d},\qquad
d_S(\omega_{i,N},\omega_{j,N})
\ge 2c_0N^{-1/d}
\quad\text{for }i\ne j.
$$
\end{lemma}
\begin{proof}
This follows from
\cite[Theorem~2]{GiganteLeopardi2017}.
\end{proof}

Recall that a probability measure on $S$ is a nonnegative Borel measure $\mu$ with $\mu(S)=1$.
For probability measures $\mu,\nu$ on $S$, define
$$
D_{\mathrm{ball}}[\mu,\nu]
=\sup_{\substack{x\in\mathbb R^n\\r>0}}
|\mu(B(x,r))-\nu(B(x,r))|.
$$
For a signed measure $\eta$ on $S$, define
$$
E_\beta(\eta)
=\iint_{S\times S}|x-y|^{-\beta}
\,d\eta(x)\,d\eta(y)
$$
whenever the integral is well-defined.

\begin{lemma}\label{lem:riesz-discrepancy}
Let $0<\beta<\gamma$, and let $\mu,\nu$ be probability measures on $S$. Suppose that for some $C_0>0$,
$$
(\mu+\nu)(B(x,r))\le C_0r^\gamma
$$
for every $x\in\mathbb R^n$ and $r>0$. Then
$$
E_\beta(\mu-\nu)
\le CD_{\mathrm{ball}}[\mu,\nu]^{1-\beta/\gamma},
$$
where $C$ depends only on $\beta$, $\gamma$, and $C_0$.
\end{lemma}

\begin{proof}
This follows from \cite[Proposition~13]{Gotz2003}.
\end{proof}

\begin{lemma}\label{lem:riesz-sobolev}
Let $0<\beta<d$. There exists $C>0$,
depending only on $\beta$ and $d$, such that
for every signed measure $\mu$ on $S$ with
$d\mu=h\,d\sigma$, where $h\in L^2(S)$ is real-valued,
$$
\|\mu\|_{H^{(\beta-d)/2}(S)}^2
\le CE_\beta(\mu).
$$
\end{lemma}

\begin{proof}
This follows from
\cite[Lemma~2.4]{MarzoMas2021}.
\end{proof}

\section{Proof of Theorem~\ref{thm:main-overshoot}}\label{sec:sieve}

We construct sets $K_N$ and measures $\mu_N$
such that each $K_N$ is a finite union of
pairwise disjoint closed balls and
\[
    \operatorname{supp}\mu_N\subset K_N\cap S,
    \qquad
    \mu_N\to\sigma
    \quad\text{strongly in }H^{-(m-1/2)}(S).
\]
We then show that
\[
    U_{K_N}(0)\longrightarrow W(0)>2.
\]

\subsection{Approximation of the surface measure}\label{sec:surface-approx}
Set
$$d=n-1,\qquad s=m-\frac{1}{2},\qquad \beta=d-2s=n-2m>0.$$
Let $\delta_N=N^{-1/d}$, and choose the sets $Q_{j,N}$ and points $\omega_{j,N}$ given by \Cref{lem:spherical-partition}, so that
\begin{equation}\label{eq:s3-2}
\operatorname{diam}_{d_S}(Q_{j,N})\le 2C_0\delta_N,\qquad d_S(\omega_{i,N},\omega_{j,N})\ge 2c_0\delta_N\quad\text{for }i\ne j.
\end{equation}
Fix
$$\frac{1}{\beta+2}<\alpha<\frac{1}{\beta},\qquad a_N=N^{-\alpha}.$$
Since $m\ge2$ and $d=n-1$, we have $\beta+2=n-2m+2<d$, and hence $\alpha>1/d$. It follows that
\begin{equation}\label{eq:s3-4}
\frac{a_N}{\delta_N}\longrightarrow0,\qquad Na_N^\beta\longrightarrow\infty,\qquad Na_N^{\beta+2}\longrightarrow0.
\end{equation}

Define
$$K_N=\bigcup_{j=1}^N\overline{B(\omega_{j,N},a_N)},$$
and write $U_N=U_{K_N}$. By \eqref{eq:s3-2} and $d_S(x,y)\le\frac{\pi}{2}|x-y|$, we have
$$|\omega_{i,N}-\omega_{j,N}|\ge\frac{4c_0}{\pi}\delta_N\qquad\text{for }i\ne j.$$
In view of $a_N/\delta_N\to0$, the closed balls $\overline{B(\omega_{j,N},2a_N)}$, $1\le j\le N$, are pairwise disjoint for all sufficiently large $N$.

Fix $0<c<1/4$, and set
$$E_{j,N}=B_S(\omega_{j,N},ca_N),\qquad \nu_{j,N}=\frac{\sigma|_{E_{j,N}}}{\sigma(E_{j,N})},\qquad \mu_N=\frac{1}{N}\sum_{j=1}^N\nu_{j,N}.$$
For any $x\in\overline{E_{j,N}}$, we have $|x-\omega_{j,N}|\le d_S(x,\omega_{j,N})\le ca_N<a_N$, so that $\overline{E_{j,N}}\subset B(\omega_{j,N},a_N)$. In particular, the sets $E_{j,N}$ are pairwise disjoint for all sufficiently large $N$.

\begin{lemma}\label{lem:surface-approx}
Let $\mu_N$ be as defined above. Then we have
\begin{equation}\label{eq:s3-7}
\mu_N\longrightarrow\sigma\qquad\text{in }H^{-s}(S).
\end{equation}
\end{lemma}

\begin{proof}
Choose $\beta<\gamma<1/\alpha<d$. We first show that
\begin{equation}\label{eq:s3-9}
(\mu_N+\sigma)(B(x,r))\le Cr^\gamma
\end{equation}
for every $x\in\mathbb R^n$ and $r>0$, uniformly for all sufficiently large $N$.
Since $E_{j,N}\subset B(\omega_{j,N},a_N)$ and $|\omega_{i,N}-\omega_{j,N}|>4a_N$ for $i\ne j$, any ball $B(x,r)$ with $0<r\le a_N$ meets at most one of the sets $E_{j,N}$. Using the surface-measure estimates
$$\sigma(E_{j,N})\ge C^{-1}a_N^d,\qquad \sigma(B(x,r))\le Cr^d,$$
we obtain, for $0<r\le a_N$,
$$\mu_N(B(x,r))\le\frac{C}{N}\left(\frac{r}{a_N}\right)^d\le\frac{C}{Na_N^\gamma}r^\gamma\le Cr^\gamma.$$

For $a_N<r\le 1$, if $B(x,r)\cap E_{j,N}\ne\varnothing$, then \eqref{eq:s3-2} and $E_{j,N}\subset B(\omega_{j,N},a_N)$ give
$$Q_{j,N}\subset B(x,r+a_N+2C_0\delta_N).$$
Since the sets $Q_{j,N}$ are pairwise disjoint
with $\sigma(Q_{j,N})=1/N$, we obtain,
for all sufficiently large $N$ and $a_N<r\le1$,
$$
\mu_N(B(x,r))
\le \sigma(B(x,r+a_N+2C_0\delta_N))
\le C(r+\delta_N)^d
\le C\left(r^d+\frac1N\right)
\le Cr^\gamma,
$$
where we used $N\delta_N^d=1$,
$\gamma<d$, and $Na_N^\gamma\to\infty$.
This, together with $\sigma(B(x,r))\le Cr^d\le Cr^\gamma$, gives $(\mu_N+\sigma)(B(x,r))\le Cr^\gamma$ for $0<r\le 1$. For $r\ge 1$, the same estimate follows from $(\mu_N+\sigma)(B(x,r))\le 2$. This proves \eqref{eq:s3-9}.

Next, let $A=B_S(\xi,\theta)$ be a spherical cap. For $0<t\le 1$, we have
$$\sigma\bigl(\{\omega\in S:d_S(\omega,\partial A)\le t\}\bigr)\le Ct.$$
Note that $d_S(\omega,y)\le 2C_0\delta_N+ca_N$ whenever $\omega\in Q_{j,N}$ and $y\in E_{j,N}$. Thus $\mathbf1_A(y)\ne\mathbf1_A(\omega)$ implies $d_S(\omega,\partial A)\le 2C_0\delta_N+ca_N$.
Since $\bigcup_{j=1}^N Q_{j,N}=S$, $\sigma(Q_{j,N})=1/N$, and each $\nu_{j,N}$ is a probability measure supported in $E_{j,N}$, we obtain
$$\begin{aligned}
|\mu_N(A)-\sigma(A)|
&=\left|\sum_{j=1}^N\int_{Q_{j,N}}\int_{E_{j,N}}\bigl(\mathbf1_A(y)-\mathbf1_A(\omega)\bigr)\,d\nu_{j,N}(y)\,d\sigma(\omega)\right|\\
&\le\sigma\bigl(\{\omega\in S:d_S(\omega,\partial A)\le 2C_0\delta_N+ca_N\}\bigr)\le C\delta_N,
\end{aligned}$$
where we used $a_N=o(\delta_N)$ in the last step. Hence
\begin{equation}\label{eq:s3-10}
\sup_A|\mu_N(A)-\sigma(A)|\le C\delta_N.
\end{equation}

Since the intersection of any Euclidean ball with $S$ is a spherical cap, the empty set, $S$, or $S$ with one point removed, \eqref{eq:s3-10} implies $D_{\mathrm{ball}}[\mu_N,\sigma]\le C\delta_N\to0$. By \Cref{lem:riesz-discrepancy} and \eqref{eq:s3-9}--\eqref{eq:s3-10}, we have
$$E_\beta(\mu_N-\sigma)\le CD_{\mathrm{ball}}[\mu_N,\sigma]^{1-\beta/\gamma}\longrightarrow0.$$

Finally, we have
$$
d(\mu_N-\sigma)=h_N\,d\sigma,\qquad
h_N=\frac{1}{N}\sum_{j=1}^N
\frac{\mathbf1_{E_{j,N}}}{\sigma(E_{j,N})}-1
\in L^2(S).
$$
Since $(\beta-d)/2=-s$, \Cref{lem:riesz-sobolev} gives
$$
\|\mu_N-\sigma\|_{H^{-s}(S)}^2
\le CE_\beta(\mu_N-\sigma)
\longrightarrow0,
$$
which proves \eqref{eq:s3-7}.
\end{proof}

\begin{lemma}\label{thm:trace-forcing}
Let $v_N\in\mathcal A_{K_N}$ satisfy
$v_N\rightharpoonup v$ in
$\mathring W^{m,2}(\mathbb R^n)$. Then
$$
\operatorname{Tr}_S v=1
\qquad\text{and}\qquad
\operatorname{Tr}_S v_N\longrightarrow1
\quad\text{in }L^2(S).
$$
\end{lemma}

\begin{proof}

Set $g_N=\operatorname{Tr}_S v_N-1$. Since $v_N\rightharpoonup v$ in $\mathring W^{m,2}(\mathbb R^n)$ and the trace operator is a bounded linear operator,
$\{g_N\}$ is bounded in $H^s(S)$. It follows from $E_{j,N}\subset B(\omega_{j,N},a_N)$ and $v_N\in\mathcal A_{K_N}$ that $g_N=0$ $\mu_N$-almost everywhere. Thus, for every $\varphi\in C^\infty(S)$,
$$\left|\int_S g_N\varphi\,d\sigma\right|=|\langle\sigma-\mu_N,g_N\varphi\rangle|\le C_\varphi\|g_N\|_{H^s(S)}\|\sigma-\mu_N\|_{H^{-s}(S)}\longrightarrow0$$
by \Cref{lem:surface-approx}. On the other hand, $g_N\rightharpoonup\operatorname{Tr}_S v-1$ in $H^s(S)$, and hence $\operatorname{Tr}_S v=1$. Since $s>0$ and $S$ is compact, the embedding $H^s(S)\Subset L^2(S)$ is compact, and hence $g_N\to0$ in $L^2(S)$.
\end{proof}

\subsection{Convergence of capacitary potentials}\label{sec:cap-pot-convergence}

\begin{lemma}\label{thm:sieve-convergence}
Let $U_N=U_{K_N}$, where $K_N$ is defined in \Cref{sec:surface-approx}. Then we have
\begin{equation}\label{eq:s3-17}
U_N\longrightarrow W\qquad\text{in }\mathring W^{m,2}(\mathbb R^n).
\end{equation}
Moreover, $U_N\to W$ uniformly on compact subsets of $B_{1/2}$.
\end{lemma}

\begin{proof}
By the regularity of $W$ established in \Cref{sec:membrane}, $D^\alpha W$ is bounded in a fixed neighborhood of $S$ for all $|\alpha|\le m$. Since $W=1$ on $S$, the mean value theorem gives
\begin{equation}\label{eq:s3-15}
|1-W(x)|\le C\operatorname{dist}(x,S)\qquad\text{for }x\text{ near }S.
\end{equation}
For sufficiently large $N$, choose $\eta_{j,N}\in C_0^\infty(B(\omega_{j,N},2a_N))$ such that $\eta_{j,N}=1$ in $B(\omega_{j,N},3a_N/2)$ and $|D^\alpha\eta_{j,N}|\le Ca_N^{-|\alpha|}$ for $|\alpha|\le m$. Since the balls $B(\omega_{j,N},2a_N)$ are pairwise disjoint, so are the supports of $\eta_{j,N}$. Set
$$W_N=W+\sum_{j=1}^N\eta_{j,N}(1-W).$$
Since $W_N=1$ in a neighborhood of $K_N$, \Cref{lem:unit-neighborhood} gives $W_N\in\mathcal A_{K_N}$. For sufficiently large $N$, the supports of $\eta_{j,N}$ lie in the neighborhood where \eqref{eq:s3-15} holds. By \eqref{eq:s3-15}, we obtain
$$\|\eta_{j,N}(1-W)\|_{\mathring W^{m,2}}^2\le Ca_N^{n-2m+2}=Ca_N^{\beta+2}.$$
Hence, by \eqref{eq:s3-4},
$$\|W_N-W\|_{\mathring W^{m,2}}^2=\sum_{j=1}^N\|\eta_{j,N}(1-W)\|_{\mathring W^{m,2}}^2\le CNa_N^{\beta+2}\longrightarrow0.$$

On the other hand, since $U_N$ is the minimizer in $\mathcal A_{K_N}$ and $W_N\in\mathcal A_{K_N}$, we have $\mathcal E_m(U_N,U_N)\le\mathcal E_m(W_N,W_N)$. As $W_N\to W$ in $\mathring W^{m,2}(\mathbb R^n)$, this implies
\begin{equation}\label{eq:s3-14}
\limsup_{N\to\infty}\mathcal E_m(U_N,U_N)\le\mathcal E_m(W,W).
\end{equation}
In particular, the sequence  $\{U_{N}\}$ is bounded in $\mathring W^{m,2}(\mathbb R^n)$.
For any subsequence  $\{U_{N_k}\}$, we may choose
a further subsequence, still denoted by
 $\{U_{N_k}\}$, such that
$U_{N_k}\rightharpoonup U$
in $\mathring W^{m,2}(\mathbb R^n)$
for some $U\in\mathring W^{m,2}(\mathbb R^n)$.
Since $U_{N_k}\in\mathcal A_{K_{N_k}}$, \Cref{thm:trace-forcing} gives $\operatorname{Tr}_S U=1$, and hence $U\in\mathcal M_S$. By the minimality of $W$ in $\mathcal M_S$, the weak lower semicontinuity of the $\mathring W^{m,2}$ norm, and \eqref{eq:s3-14}, we obtain
$$\mathcal E_m(W,W)\le\mathcal E_m(U,U)\le\liminf_{k\to\infty}\mathcal E_m(U_{N_k},U_{N_k})\le\limsup_{k\to\infty}\mathcal E_m(U_{N_k},U_{N_k})\le\mathcal E_m(W,W).$$
By uniqueness of the minimizer in $\mathcal M_S$, we have $U=W$ and $\|U_{N_k}\|_{\mathring W^{m,2}}\to\|W\|_{\mathring W^{m,2}}$, and hence $U_{N_k}\to W$ in $\mathring W^{m,2}(\mathbb R^n)$. Since the subsequence was arbitrary, \eqref{eq:s3-17} follows.

Finally, since $a_N\to0$, we have $K_N\cap B_{1/2}=\varnothing$ for sufficiently large $N$, and hence $U_N$ and $W$ are $m$-harmonic in $B_{1/2}$. For any compact subset $E\subset B_{1/2}$, the interior estimates~\cite[Section~2.5]{GazzolaGrunauSweers2010} and the Sobolev inequality, together with \eqref{eq:s3-17}, yield
$$\sup_{x\in E}|U_N(x)-W(x)|\le C\|U_N-W\|_{\mathring W^{m,2}}\longrightarrow0.$$
\end{proof}

\begin{proof}[\textnormal{\textbf{Proof of \Cref{thm:main-overshoot}}}]
By \eqref{eq:s2-24}, we have $W(0)>2$.
\Cref{thm:sieve-convergence} gives
$U_{K_N}(0)\to W(0)$.
Thus $U_{K_N}(0)>2$ for all sufficiently large $N$.
Taking $K=K_N$ for such an $N$, we complete the proof.
\end{proof}

\begin{remark}\label{rem:uniform-sharpness}
The constant $2$ in Maz'ya's estimate $U_F<2$  cannot be replaced by a smaller constant independent of $m$ and $n$. Indeed, for $m\ge 3$ and $n=2m+2$, it
follows from \eqref{eq:s2-23} that $W(0)=2-1/m$. By the construction in Section~\ref{sec:surface-approx}
and Lemma~\ref{thm:sieve-convergence}, $U_{K_N}(0)\to W(0)$ as $N\to\infty$.
Letting $m\to\infty$ yields the desired conclusion.
\end{remark}

\section{Proofs of Theorems~\ref{thm:main-unbounded} and~\ref{thm:main-green}}\label{sec:amplifier}
In this section, we first establish an asymptotic relation
between the Dirichlet solutions in $G$
and $G\setminus rK$ as $r\downarrow0$.
We then use this relation to construct an unbounded solution. The same construction also yields
counterexamples concerning Green functions.

Let $G\subset\mathbb R^n$ be a bounded open set,
$E\Subset G$ a compact set, and $h\in H_0^m(G)$.
Set
\[
\mathcal A_G(E,h)
=h+H_0^m(G\setminus E)\subset H_0^m(G).
\]
By the Hilbert projection theorem, the minimization problem
$$\inf_{v\in\mathcal A_G(E,h)}\mathcal E_m(v,v)$$
has a unique minimizer $W_{E,h}^G\in\mathcal A_G(E,h)$, characterized by
\begin{equation}
\mathcal E_m(W_{E,h}^G,\varphi)=0
\qquad\text{for every }
\varphi\in H_0^m(G\setminus E).
\label{eq:s4-2}
\end{equation}

\subsection{Asymptotic behavior of solutions}\label{sec:solution-asymptotics}

\begin{lemma}\label{lem:expanding-mosco}
Let $G\subset\mathbb R^n$ be a bounded open set
with $0\in G$, and let $K\subset\mathbb R^n$
be a compact set. Set
\[
G_r=r^{-1}G,\qquad V_r=H_0^m(G_r\setminus K),
\]
with zero extension to $\mathbb R^n$.
As $r\downarrow0$, the spaces $V_r$ converge to
$\mathring W^{m,2}(\mathbb R^n\setminus K)$
in the sense of Mosco in
$\mathring W^{m,2}(\mathbb R^n)$.
\end{lemma}

\begin{proof}
Set $V=\mathring W^{m,2}(\mathbb R^n\setminus K)$.
Since $0\in G$, there exists $\rho>0$
such that $B_\rho\subset G$, and hence
$B_{\rho/r}\subset G_r$ for $r>0$.

For ($\mathrm{M1}_r$), let $v\in V$
and set $v_r=P_{V_r}v$,
where $P_{V_r}$ denotes the orthogonal projection
onto $V_r$ in $\mathring W^{m,2}(\mathbb R^n)$.
For any $\varphi\in C_0^\infty(\mathbb R^n\setminus K)$,
we have
$\varphi\in H_0^m(B_{\rho/r}\setminus K)\subset V_r$
for all sufficiently small $r>0$, so that
\[
\limsup_{r\downarrow0}
\|v-v_r\|_{\mathring W^{m,2}}
\le\|v-\varphi\|_{\mathring W^{m,2}}.
\]
By a density argument, we obtain $v_r\to v$
in $\mathring W^{m,2}(\mathbb R^n)$,
which gives ($\mathrm{M1}_r$).

For ($\mathrm{M2}_r$), let $r_j\to0^+$
and $v_j\in V_{r_j}$ be such that
$v_j\rightharpoonup v$
in $\mathring W^{m,2}(\mathbb R^n)$.
Since $V_{r_j}\subset V$ and $V$ is weakly closed,
we obtain $v\in V$.
Thus ($\mathrm{M2}_r$) is satisfied.
\end{proof}

\begin{lemma}\label{prop:local-blowup}
Let $G\subset\mathbb R^n$ be a bounded open set
with $0\in G$, and $K\subset\mathbb R^n$
a compact set. Suppose that $v\in H_0^m(G)$
is smooth near $0$, with $A=v(0)\ne0$. For sufficiently small $r>0$, set $G_r=r^{-1}G$, $H_r=W_{rK,v}^G$, and
$$\widetilde H_r(y)=\frac{H_r(ry)}{A},$$
where $\widetilde H_r$ is extended by zero outside $G_r$. Then, as $r\downarrow0$,
\begin{align}
\widetilde H_r &\longrightarrow U_K \qquad\text{in }\mathring W^{m,2}(\mathbb R^n),\label{eq:s4-4}\\
\mathcal E_m(H_r,H_r) &= A^2r^{n-2m}\bigl(\operatorname{cap}_m(K)+o(1)\bigr).\label{eq:s4-5}
\end{align}
Moreover, $\widetilde H_r\to U_K$ uniformly on compact subsets of $\mathbb R^n\setminus K$. In particular, if $0\notin K$, then
\begin{equation}\label{eq:s4-5a}
\frac{H_r(0)}{A}\longrightarrow U_K(0).
\end{equation}
\end{lemma}

\begin{proof}
Let $V_r=H_0^m(G_r\setminus K)$ as in \Cref{lem:expanding-mosco}.
Choose $\eta\in\mathcal C_K$.
Since $0\in G$, we have
$\operatorname{supp}\eta\subset G_r$
for all sufficiently small $r>0$. Set
$$v_r(y)=\frac{v(ry)}{A},\qquad b_r=\eta v_r.$$
Since $v$ is smooth near $0$ and $v(0)=A$, Taylor's formula and Leibniz' rule give
\begin{equation}\label{eq:br-convergence}
b_r\longrightarrow\eta\qquad\text{in }\mathring W^{m,2}(\mathbb R^n).
\end{equation}
Moreover, $v_r\in H_0^m(G_r)$, and since $\eta\in\mathcal C_K$,
$$b_r-v_r=(\eta-1)v_r\in H_0^m(G_r\setminus K)=V_r.$$
By \eqref{eq:s4-2} and rescaling, $\widetilde H_r$
is the minimizer of $\mathcal E_m(\cdot,\cdot)$
in $b_r+V_r$. \Cref{lem:expanding-mosco,lem:mosco-minimizers},
together with \eqref{eq:br-convergence}, give \eqref{eq:s4-4}.

Using \eqref{eq:s4-4} and the definition of $\widetilde H_r$, we obtain
$$\mathcal E_m(H_r,H_r)=A^2r^{n-2m}\mathcal E_m(\widetilde H_r,\widetilde H_r)=A^2r^{n-2m}\bigl(\operatorname{cap}_m(K)+o(1)\bigr)$$
as $r\downarrow0$, which gives \eqref{eq:s4-5}.

Finally, let $E$ be a compact subset of
$\mathbb R^n\setminus K$.
For sufficiently small $r>0$,
the function $\widetilde H_r-U_K$ is $m$-harmonic
in a fixed neighborhood of $E$.
The interior estimates and the Sobolev inequality,
together with \eqref{eq:s4-4}, give
$$
\sup_{x\in E}|\widetilde H_r(x)-U_K(x)|
\le C\|\widetilde H_r-U_K\|_{\mathring W^{m,2}(\mathbb R^n)}
\longrightarrow0.
$$
\end{proof}

\subsection{Construction of an unbounded solution}\label{sec:construction-unbounded}

Let $u\in\mathring W^{m,2}(\mathbb R^n)$ satisfy
$(-\Delta)^m u=0$ in $B(x,2r)$.
The interior estimates for $m$-harmonic functions
and the homogeneous Sobolev inequality give
\begin{equation}
|u(x)|\le Cr^{m-n/2}
\|u\|_{\mathring W^{m,2}(\mathbb R^n)},
\label{eq:amplifier:20b}
\end{equation}
where $C=C(n,m)$.

\begin{lemma}\label{thm:amplifier}
Let $K\subset\mathbb R^n\setminus\{0\}$ be a compact set, and set $M=U_K(0)$. Assume $|1-M|>1$. Then there exist a bounded connected domain $\Omega^*\subset\mathbb R^n$ and $f\in C_0^\infty(\Omega^*)$, with $f\ge0$ and $f\not\equiv0$, such that the variational solution $u\in H_0^m(\Omega^*)$ of $(-\Delta)^m u=f$ does not belong to $L^\infty(\Omega^*)$. More precisely, if $M>2$, then
$$\operatorname*{ess\,sup}_{\Omega^*}u=+\infty,\qquad \operatorname*{ess\,inf}_{\Omega^*}u=-\infty,$$
while if $M<0$, then $\operatorname*{ess\,sup}_{\Omega^*}u=+\infty$.
\end{lemma}

\begin{proof}
Let $a=\operatorname{dist}(0,K)$, $b=\max_{x\in K}|x|$, and $\varsigma=\operatorname{sgn}(1-M)$. Since $0\notin K$, $U_K$ is $m$-harmonic and hence smooth near $0$. The inequality $\varsigma(1-U_K(0))=|1-M|>1$ and continuity give a point $y_*\in B_{a/8}\setminus\{0\}$ such that $\varsigma(1-U_K(y_*))>1$.

Choose $1<\rho<|1-M|$ and $1<\kappa<\varsigma(1-U_K(y_*))$, and take $c>0$ sufficiently small so that $B(y_*,4c)\cap K=\varnothing$ and $4c<|y_*|$. Fix $f\in C_0^\infty(B_3\setminus\overline{B_2})$ with $f\ge0$ and $f\not\equiv0$, and let $\psi\in H_0^m(B_4)$ be the variational solution of $(-\Delta)^m\psi=f$ in $B_4$. Then $\psi$ is smooth in $B_4$, and the Green representation formula together with the positivity of the Dirichlet Green function for $(-\Delta)^m$ in a ball~\cite{Boggio1905} gives
$$\psi(0)=\int_{B_4}G_{B_4}(0,y)f(y)\,dy>0.$$

We construct the radii $r_j>0$ inductively, requiring $br_1<1$ and, for $j\ge2$,
\begin{equation}\label{eq:s4-7}
br_j<\min\left\{\frac{a}{4}r_{j-1},\,\frac{|y_*|}{8}r_{j-1}\right\}.
\end{equation}
For such radii, set $\Omega_0=B_4$, $H_0=H_0^m(\Omega_0)$, and
$$K_j=r_jK,\qquad \Omega_j=B_4\setminus\bigcup_{k=1}^j K_k,\qquad H_j=H_0^m(\Omega_j)\quad\text{for }j\ge1.$$
By zero extension, each $H_j$ is a closed subspace of $\mathring W^{m,2}(\mathbb R^n)$.

Let $u_j=P_{H_j}\psi$. Then
\begin{equation}\label{eq:s4-8}
\mathcal E_m(u_j,\phi)=\langle f,\phi\rangle\qquad\text{for every }\phi\in H_j.
\end{equation}
In particular, $u_0=\psi$ with $u_0(0)>0$. Moreover, $u_j$ is $m$-harmonic in $\Omega_j\setminus\operatorname{supp}f$ and hence smooth near $0$.

Define
$$A_j=u_j(0)\quad\text{for }j\ge0,\qquad h_j=u_{j-1}-u_j\quad\text{for }j\ge1.$$
For $j\ge2$, \eqref{eq:s4-7} gives $K_j\subset B_{ar_{j-1}/4}$. Since $K_1,\dots,K_{j-1}$ lie outside $B_{ar_{j-1}}$, we have $K_j\Subset\Omega_{j-1}$; for $j=1$ this follows from $br_1<1$.
Hence $H_j=H_0^m(\Omega_{j-1}\setminus K_j)
\subset H_{j-1}$, so that
$h_j\in\mathcal A_{\Omega_{j-1}}(K_j,u_{j-1})$.
Moreover, \eqref{eq:s4-8} gives
\begin{equation}\label{eq:s4-9}
\mathcal E_m(h_j,\phi)=0\qquad\text{for every }\phi\in H_j.
\end{equation}
This, together with \eqref{eq:s4-2}, yields $h_j=W_{K_j,u_{j-1}}^{\Omega_{j-1}}$.

Suppose that $r_1,\dots,r_{j-1}$ have been chosen with $A_{j-1}\ne0$. Applying \Cref{prop:local-blowup} to $u_{j-1}$ on $\Omega_{j-1}$, we obtain
$h_j(0)/A_{j-1}\to M$ as $r_j\downarrow0$.
Since $A_j=A_{j-1}-h_j(0)$, it follows that
\begin{equation}\label{eq:s4-10}
\frac{A_j}{A_{j-1}}=1-\frac{h_j(0)}{A_{j-1}}\longrightarrow 1-M.
\end{equation}
Set $x_j=r_jy_*$.
Since $u_{j-1}(x_j)\to A_{j-1}$
as $r_j\downarrow0$, we obtain
\begin{equation}\label{eq:s4-11}
\frac{u_j(x_j)}{A_{j-1}}
=\frac{u_{j-1}(x_j)}{A_{j-1}}
-\frac{h_j(x_j)}{A_{j-1}}
\longrightarrow 1-U_K(y_*).
\end{equation}
Moreover, by \eqref{eq:s4-5},
\begin{equation}\label{eq:s4-12}
\|h_j\|_{\mathring W^{m,2}}^2=A_{j-1}^2r_j^{n-2m}\bigl(\operatorname{cap}_m(K)+o(1)\bigr).
\end{equation}
Since $n>2m$, \eqref{eq:s4-12} implies $\|h_j\|_{\mathring W^{m,2}}\to0$ as $r_j\downarrow0$. In view of \eqref{eq:s4-10}--\eqref{eq:s4-11}, there exists $r_j>0$ sufficiently small, satisfying $br_1<1$ if $j=1$ or \eqref{eq:s4-7} if $j\ge2$, such that
\begin{align}
\varsigma\frac{A_j}{A_{j-1}} &\ge\rho,\label{eq:s4-13}\\
\varsigma\frac{u_j(x_j)}{A_{j-1}} &\ge\kappa,\label{eq:s4-14}\\
Cr_i^{m-n/2}\|h_j\|_{\mathring W^{m,2}} &\le2^{-j}|A_{i-1}|\qquad\text{for }1\le i<j.\label{eq:s4-15}
\end{align}
In particular, \eqref{eq:s4-13} and $A_0>0$ ensure that $A_j\ne0$ and
\begin{equation}\label{eq:s4-16}
|A_j|\ge A_0\rho^j,\qquad \operatorname{sgn}(A_j)=\varsigma^j.
\end{equation}

Define
$$
F_\infty=\{0\}\cup\bigcup_{j=1}^\infty K_j,
\qquad
\Omega=B_4\setminus F_\infty.
$$
Since $r_j\to0$ and $n>2m$,
a standard cutoff argument gives
$
H_0^m(\Omega)=\bigcap_{j=1}^\infty H_j.
$

Since $u_j=P_{H_j}\psi$ and the closed subspaces
$H_j$ are decreasing, we obtain
\begin{equation}\label{eq:s4-17}
u_j\longrightarrow u=P_{H_0^m(\Omega)}\psi
\qquad\text{in }\mathring W^{m,2}(\mathbb R^n).
\end{equation}

By \eqref{eq:s4-8} and \eqref{eq:s4-17}, we obtain
\begin{equation}\label{eq:s4-18}
\mathcal E_m(u,\phi)=\langle f,\phi\rangle
\qquad\text{for every }\phi\in H_0^m(\Omega).
\end{equation}

We next estimate $u$ at $x_i=r_iy_*$. The choices of $y_*$ and $c$ imply
$$\overline{B(x_i,2cr_i)}\cap K_i=\varnothing,\qquad \overline{B(x_i,2cr_i)}\subset B_{br_i}\setminus\overline{B_{|x_i|/2}}.$$
By \eqref{eq:s4-7}, $K_k\subset B_{|x_i|/8}$ for $k>i$, while $K_k\cap B_{br_i}=\varnothing$ for $k<i$. This, together with $br_i<1$, implies $B(x_i,2cr_i)\Subset\Omega$. Thus $C_0^\infty(B(x_i,2cr_i))\subset H_j$ for every $j\ge1$, and \eqref{eq:s4-9} shows that $h_j$ is $m$-harmonic in $B(x_i,2cr_i)$.
For $j>i$, \eqref{eq:amplifier:20b} and \eqref{eq:s4-15} give
\begin{equation}\label{eq:s4-19}
|h_j(x_i)|\le Cr_i^{m-n/2}\|h_j\|_{\mathring W^{m,2}}\le 2^{-j}|A_{i-1}|.
\end{equation}
For each fixed $i$, \eqref{eq:s4-8} and
\eqref{eq:s4-18} show that $u_N-u$ is
$m$-harmonic in $B(x_i,2cr_i)$.
The interior estimates, together with
\eqref{eq:s4-17}, give
$u_N(x_i)\to u(x_i)$ as $N\to\infty$. Since $u_N(x_i)=u_i(x_i)-\sum_{k=i+1}^N h_k(x_i)$ for $N>i$, letting $N\to\infty$ and using \eqref{eq:s4-19} yields
$$|u(x_i)-u_i(x_i)|\le\sum_{k=i+1}^\infty 2^{-k}|A_{i-1}|=2^{-i}|A_{i-1}|.$$
In view of \eqref{eq:s4-14} and \eqref{eq:s4-16}, it follows that
\begin{equation}\label{eq:s4-20}
\varsigma^i u(x_i)\ge(\kappa-2^{-i})|A_{i-1}|\ge A_0(\kappa-2^{-i})\rho^{i-1}\longrightarrow+\infty\qquad\text{as }i\to\infty.
\end{equation}

Finally, since $K_j\subset B_1$ for all $j$,
the connected set $B_3\setminus\overline{B_2}$
lies in a connected component $\Omega^*$ of
$\Omega$. Hence,
$\operatorname{supp}f\subset\Omega^*$.
By an approximation argument, we have
$u|_{\Omega^*}\in H_0^m(\Omega^*)$.
Similarly, for any other connected component
$C$ of $\Omega$, the zero extension $u_C$
of $u|_C$ to $\Omega$ belongs to
$H_0^m(\Omega)$. Since $f=0$ on $C$, \eqref{eq:s4-18} gives
$$\mathcal E_m(u_C,u_C)=\mathcal E_m(u,u_C)=\langle f,u_C\rangle=0.$$

Thus $u=0$ on $\Omega\setminus\Omega^*$.
By \eqref{eq:s4-18}, $u|_{\Omega^*}$ is the
variational solution of $(-\Delta)^m u=f$
in $\Omega^*$. Since $\kappa>1$ and $u$ is smooth near each $x_i$,
it follows from \eqref{eq:s4-20} that
$x_i\in\Omega^*$ for every $i\ge1$. Hence $u|_{\Omega^*}$ is unbounded above on $\Omega^*$, and also unbounded below if $\varsigma=-1$.

\end{proof}

\begin{proof}[\textnormal{\textbf{Proof of \Cref{thm:main-unbounded}}}]
By Theorem~\ref{thm:main-overshoot}, there exists a finite union
$K\Subset\mathbb{R}^n\setminus\{0\}$ of pairwise disjoint closed balls such
that $U_K(0)>2$. Applying \Cref{thm:amplifier} to this $K$ and
using the notation in its proof, we have $\varsigma=-1$, and the variational
solution on $\Omega^*$ is unbounded both above and below.

For every $i\ge1$, the point $x_i=r_iy_*$ belongs to
$\Omega^*$. Since $B(x_i,2cr_i)$ is a connected subset
of $\Omega$ containing $x_i$, it is also contained in
$\Omega^*$. Since $x_i\to0$ and $0\notin\Omega^*$,
we have $0\in\partial\Omega^*$. Moreover,
\[
\operatorname{dist}(x_i,\partial\Omega^*)
\ge 2cr_i
=\frac{2c}{|y_*|}|x_i|.
\]
Taking $e=y_*/|y_*|$ and $t_i=r_i|y_*|$, we obtain the nontangential sequence
from \eqref{eq:s4-20}.

It remains to verify \eqref{eq:wiener} for $\Omega^*$. Since $b=\max_{x\in K}|x|$ and $K_j=r_jK$, we have
$K_j\subset\overline{B_r}\setminus\Omega^*$
for $br_j<r<2br_j$. By monotonicity and rescaling,
\[
\frac{\operatorname{cap}_m(\overline{B_r}\setminus\Omega^*)}
{\operatorname{cap}_m(\overline{B_r})}
\ge
\frac{r_j^{n-2m}\operatorname{cap}_m(K)}
{r^{n-2m}\operatorname{cap}_m(\overline{B_1})}
\ge
\frac{\operatorname{cap}_m(K)}
{(2b)^{n-2m}\operatorname{cap}_m(\overline{B_1})}>0,
\]
where $\operatorname{cap}_m(K)>0$ since $U_K(0)>2$. By \eqref{eq:s4-7}, there exists $j_0\ge 1$ such that
the intervals $(br_j,2br_j)$ ($j\ge j_0$) are pairwise disjoint and contained in $(0,1)$.
Consequently,
\[
\begin{aligned}
\int_0^1
\frac{\operatorname{cap}_m(\overline{B_r}\setminus\Omega^*)}
{\operatorname{cap}_m(\overline{B_r})}\,\frac{dr}{r}
&\ge
\frac{\operatorname{cap}_m(K)}
{(2b)^{n-2m}\operatorname{cap}_m(\overline{B_1})}
\sum_{j=j_0}^\infty
\int_{br_j}^{2br_j}\frac{dr}{r}\\
&=
\frac{\operatorname{cap}_m(K)}
{(2b)^{n-2m}\operatorname{cap}_m(\overline{B_1})}
\sum_{j=j_0}^\infty\log 2
=\infty.
\end{aligned}
\]
This proves \eqref{eq:wiener} and completes the proof.
\end{proof}

\subsection{Counterexamples concerning Green functions}\label{sec:green-consequences}

\begin{proof}[\textnormal{\textbf{Proof of \Cref{thm:main-green}}}]
By \Cref{thm:main-overshoot}, there exists a finite union $K\Subset\mathbb R^n\setminus\{0\}$ of pairwise disjoint closed balls such that $U_K(0)>2$. For this $K$, let $f$, $\Omega_j$, and $u_j$ be as in the proof of \Cref{thm:amplifier}. Then each $\Omega_j$ is a smooth bounded connected domain satisfying $0\in\Omega_j\subset B_4$, and \eqref{eq:s4-16} gives $(-1)^j u_j(0)\to+\infty$.

Let $G_j=G_{\Omega_j}$. By the Green representation formula,
$$u_j(0)=\int_{\Omega_j}G_j(0,y)f(y)\,dy=\int_{\operatorname{supp}f}G_j(0,y)f(y)\,dy.$$
Since $\operatorname{supp}f\subset B_3\setminus\overline{B_2}$, $G_j(0,\cdot)$ is continuous on $\operatorname{supp}f$. Thus there exist $y_j^+,y_j^-\in\operatorname{supp}f$ such that
$$G_{2j}(0,y_j^+)=\max_{y\in\operatorname{supp}f}G_{2j}(0,y),\qquad G_{2j+1}(0,y_j^-)=\min_{y\in\operatorname{supp}f}G_{2j+1}(0,y).$$
Since $f\ge0$ and $f\not\equiv0$, it follows that
$$G_{2j}(0,y_j^+)\ge\frac{u_{2j}(0)}{\|f\|_{L^1(\mathbb R^n)}}\longrightarrow+\infty,\qquad G_{2j+1}(0,y_j^-)\le\frac{u_{2j+1}(0)}{\|f\|_{L^1(\mathbb R^n)}}\longrightarrow-\infty.$$
Setting $\Omega_j^+=\Omega_{2j}$ and $\Omega_j^-=\Omega_{2j+1}$, we obtain \eqref{eq:s1-10}. Moreover, since $2\le|y_j^+|\le3$ and $n>2m$,
$$|y_j^+|^{n-2m}|G_{\Omega_j^+}(0,y_j^+)|\ge 2^{n-2m}|G_{\Omega_j^+}(0,y_j^+)|\longrightarrow+\infty,$$
which gives \eqref{eq:s1-11} and completes the proof.
\end{proof}

\section{Proofs of Theorems~\ref{thm:main-nonsubadd} and~\ref{thm:main-capacity-gap}}
\label{sec:potential-amplifier}
In this section we establish an asymptotic formula for $\operatorname{cap}_m(K\cup L_r)$ as $r\downarrow 0$, and use it to prove \Cref{thm:main-nonsubadd,thm:main-capacity-gap}.

For a compact set $E\subset\mathbb R^n$, write
$V_E=\mathring W^{m,2}(\mathbb R^n\setminus E)$.
Let $E,F\subset\mathbb R^n$ be disjoint compact sets,
and choose $\chi\in\mathcal C_E$ with
$\supp\chi\subset\mathbb R^n\setminus F$.
For any $u\in V_E\cap V_F$, there exist
$\varphi_j\in C_0^\infty(\mathbb R^n\setminus E)$
and $\psi_j\in C_0^\infty(\mathbb R^n\setminus F)$
such that $\varphi_j,\psi_j\to u$
in $\mathring W^{m,2}(\mathbb R^n)$.
Then
\[
h_j=\chi\varphi_j+(1-\chi)\psi_j
\in C_0^\infty(\mathbb R^n\setminus(E\cup F)),
\]
and $h_j\to u$ in $\mathring W^{m,2}(\mathbb R^n)$.
Thus $V_E\cap V_F\subset V_{E\cup F}$.
The reverse inclusion is immediate.

Choose $\eta\in\mathcal C_{E\cup F}$.
Since $\mathcal A_E=\eta+V_E$,
$\mathcal A_F=\eta+V_F$, and
$\mathcal A_{E\cup F}=\eta+V_{E\cup F}$, it follows that
\begin{equation}
V_{E\cup F}=V_E\cap V_F,\qquad
\mathcal A_{E\cup F}=\mathcal A_E\cap\mathcal A_F.
\label{sec6:intersection}
\end{equation}

\subsection{Failure of subadditivity}
\label{sec:proof-nonsubadd}
\begin{lemma}\label{lem:interaction}
Let $K,L\subset\mathbb R^n$ be compact sets with $\operatorname{cap}_m(L)>0$, and $x_0\in\mathbb R^n\setminus K$. Set $L_r=x_0+rL$ for $r>0$. Then, as $r\downarrow0$,
\begin{equation}\label{sec6:blowup}
U_{K\cup L_r}(x_0+r\,\cdot)-U_K(x_0+r\,\cdot)\longrightarrow(1-U_K(x_0))U_L
\end{equation}
in $\mathring W^{m,2}(\mathbb R^n)$ and uniformly on compact subsets of $\mathbb R^n\setminus L$. Moreover,
\begin{equation}\label{sec6:capacity-asymptotic}
\operatorname{cap}_m(K\cup L_r)=\operatorname{cap}_m(K)+|1-U_K(x_0)|^2\operatorname{cap}_m(L_r)+o(\operatorname{cap}_m(L_r)).
\end{equation}
\end{lemma}
\begin{proof}
By translation, we may assume that $x_0=0$. For sufficiently small $r>0$, $K\cap rL=\varnothing$. Set $z_r=U_{K\cup rL}-U_K$. By \eqref{sec6:intersection}, $z_r\in V_K$; since $U_K\perp V_K$, we have
\begin{equation}\label{sec6:energy-identity}
\operatorname{cap}_m(K\cup rL)=\operatorname{cap}_m(K)+\mathcal E_m(z_r,z_r).
\end{equation}
Moreover, since $U_{K\cup rL}$ minimizes $\mathcal E_m(u,u)$ in $\mathcal A_{K\cup rL}=\mathcal A_K\cap\mathcal A_{rL}$,
$$\mathcal E_m(z_r,z_r)=\min\{\mathcal E_m(v,v):v\in V_K,\ U_K+v\in\mathcal A_{rL}\}.$$

Choose $\eta\in\mathcal C_L$ and set
$$w_r(y)=z_r(ry),\qquad A=1-U_K(0),\qquad b_r(y)=\eta(y)(1-U_K(ry)).$$
Since $U_K$ is smooth near $0$ and
$\operatorname{dist}(0,K/r)\to\infty$, we have
\[
b_r\in V_{K/r}\quad\text{for sufficiently small }r,
\qquad
b_r\longrightarrow A\eta
\quad\text{in }\mathring W^{m,2}(\mathbb R^n).
\]
By rescaling, $w_r$ is the minimizer of $\mathcal E_m(\cdot,\cdot)$ in
$$\{v\in V_{K/r}:v-b_r\in V_L\}=b_r+(V_{K/r}\cap V_L).$$

We claim that $V_{K/r}\cap V_L$ converges to $V_L$ in the sense of Mosco as $r\downarrow0$. For $v\in V_L$, let $v_r=P_{V_{K/r}\cap V_L}v$. Every $\varphi\in C_0^\infty(\mathbb R^n\setminus L)$ belongs to $V_{K/r}\cap V_L$ for all sufficiently small $r>0$. Thus
$$\limsup_{r\downarrow0}\|v-v_r\|_{\mathring W^{m,2}}\le\|v-\varphi\|_{\mathring W^{m,2}}.$$
By the density of $C_0^\infty(\mathbb R^n\setminus L)$ in $V_L$, we obtain $v_r\to v$ in $\mathring W^{m,2}(\mathbb R^n)$, which gives ($\mathrm{M1}_r$). Since $V_{K/r}\cap V_L\subset V_L$ and $V_L$ is weakly closed, ($\mathrm{M2}_r$) is satisfied. By \Cref{lem:mosco-minimizers},
$$w_r\longrightarrow AU_L\qquad\text{in }\mathring W^{m,2}(\mathbb R^n),$$
which proves the convergence in $\mathring W^{m,2}(\mathbb R^n)$ in \eqref{sec6:blowup}.

For any $Q\Subset\mathbb R^n\setminus L$, choose a bounded open set $O$ with $Q\Subset O\Subset\mathbb R^n\setminus L$. Note that for $r>0$ sufficiently small, $w_r-AU_L$ is $m$-harmonic in $O$. It follows from the interior estimates and the Sobolev inequality that
$$\begin{aligned}
\sup_{y\in Q}|w_r(y)-AU_L(y)|
&\le C\|w_r-AU_L\|_{\mathring W^{m,2}(\mathbb R^n)}
\longrightarrow0.
\end{aligned}$$

Finally, by rescaling, we obtain
$$\begin{aligned}
\mathcal E_m(z_r,z_r)
&=r^{n-2m}\mathcal E_m(w_r,w_r)\\
&=r^{n-2m}\bigl(A^2\operatorname{cap}_m(L)+o(1)\bigr)\\
&=A^2\operatorname{cap}_m(rL)+o(\operatorname{cap}_m(rL)).
\end{aligned}$$
This, together with \eqref{sec6:energy-identity}, yields \eqref{sec6:capacity-asymptotic}.
\end{proof}

\begin{lemma}\label{lem:unbounded-potentials}
Suppose that either \eqref{eq:s1-5} or \eqref{eq:s1-6} holds. For every $M>0$, there exists a compact set $E\subset\mathbb R^n\setminus\{0\}$, which is a finite union of pairwise disjoint closed balls, such that
$$|1-U_E(0)|>M.$$
\end{lemma}

\begin{proof}
Let $K$ be the compact set from \Cref{thm:main-overshoot}, and set $M_0=U_K(0)>2$, $a=\operatorname{dist}(0,K)>0$, and $b=\max_{x\in K}|x|$. Fix $1<\rho<M_0-1$. We inductively choose radii $r_j>0$, with $r_1=1$, and set
$$F_j=\bigcup_{k=1}^j r_kK,\qquad A_j=U_{F_j}(0).$$
Thus $F_1=K$ and $A_1=M_0$.

Suppose that $r_1,\dots,r_{j-1}$ have been chosen with $A_{j-1}\ne1$. Applying \Cref{lem:interaction} to $F_{j-1}$ and $K$ at $x_0=0$, we obtain
$$\frac{U_{F_{j-1}\cup rK}(0)-1}{A_{j-1}-1}=1+\frac{U_{F_{j-1}\cup rK}(0)-A_{j-1}}{A_{j-1}-1}\longrightarrow 1-M_0<-\rho$$
as $r\downarrow0$. We may choose $r_j>0$ sufficiently small so that
$$br_j<\frac a4r_{j-1},\qquad \frac{A_j-1}{A_{j-1}-1}<-\rho.$$
It follows by induction that the sets $r_kK$ are pairwise disjoint and
$$(-1)^{j-1}(A_j-1)\ge(M_0-1)\rho^{j-1}.$$
Choosing $j$ sufficiently large
 so that $(M_0-1)\rho^{j-1}>M$ and taking $E=F_j$, we complete the proof.
\end{proof}

\begin{proof}[\textnormal{\textbf{Proof of \Cref{thm:main-nonsubadd}}}]
Let $M>1$. By \Cref{lem:unbounded-potentials}, choose
$E$ such that $|1-U_E(0)|^2>M$. Set $L=\overline{B_1}$.
Then \Cref{lem:interaction} gives
\[
\lim_{r\downarrow0}
\frac{\capm(E\cup rL)-\capm(E)}{\capm(rL)}
=|1-U_E(0)|^2>M.
\]
Since $0\notin E$, we have $E\cap rL=\varnothing$ and
$$\operatorname{cap}_m(E\cup rL)>\operatorname{cap}_m(E)+M\operatorname{cap}_m(rL)$$
for $r>0$ sufficiently small. Taking $F=rL$ completes the proof.
\end{proof}

\subsection{Comparison of capacities}\label{sec:capacity-gap}

\begin{proof}[\textnormal{\textbf{Proof of \Cref{thm:main-capacity-gap}}}]
Let $K$ be the compact set given by \Cref{thm:main-overshoot}, with $M=U_K(0)>2$, and set $L=\overline{B_1}$ and $H_r=K\cup rL$. Since $0\notin K$, $H_r$ is a finite union of pairwise disjoint closed balls for all sufficiently small $r>0$. By \eqref{eq:capacity-comparison} and the subadditivity of $\operatorname{R}_{2m}$,
$$\operatorname{Cap}_{m,\mathbb R^n}^{\ge}(H_r)\le\operatorname{R}_{2m}(H_r)\le\operatorname{R}_{2m}(K)+\operatorname{R}_{2m}(rL)\le\operatorname{cap}_m(K)+\operatorname{cap}_m(rL).$$
On the other hand, \Cref{lem:interaction} gives
$$\operatorname{cap}_m(H_r)-\operatorname{cap}_m(K)-\operatorname{cap}_m(rL)=\bigl(M(M-2)+o(1)\bigr)\operatorname{cap}_m(rL).$$ Since $M>2$, it follows that
$$\operatorname{cap}_m(H_r)-\operatorname{R}_{2m}(H_r)\ge\bigl(M(M-2)+o(1)\bigr)\operatorname{cap}_m(rL)>0$$
for $r>0$ sufficiently small. Fix such an $r$ and set $H=H_r$.
Then we obtain
$$\operatorname{Cap}_{m,\mathbb R^n}^{\ge}(H)\le\operatorname{R}_{2m}(H)<\operatorname{cap}_m(H),$$
which gives \eqref{eq:s1-13}. Finally, for $H_\varepsilon=\varepsilon H$ with $\varepsilon>0$, rescaling gives
$$\frac{\operatorname{Cap}_{m,\mathbb R^n}^{\ge}(H_\varepsilon)}{\operatorname{cap}_m(H_\varepsilon)}=\frac{\varepsilon^{n-2m}\operatorname{Cap}_{m,\mathbb R^n}^{\ge}(H)}{\varepsilon^{n-2m}\operatorname{cap}_m(H)}=\frac{\operatorname{Cap}_{m,\mathbb R^n}^{\ge}(H)}{\operatorname{cap}_m(H)}<1.$$
This yields \eqref{eq:s1-14}.
\end{proof}

\section*{Acknowledgements}
The author would like to thank Professor Jun Geng for introducing this problem and for his valuable suggestions on the proofs and the preparation of the manuscript.

The author used ChatGPT during the preparation of this work to assist with the
development and presentation of the mathematical arguments and the drafting and revision of the
text. All mathematical arguments and proofs have been carefully verified by the author, who takes full responsibility for the content of the paper.

\bigskip
\noindent
\textbf{Honggang Li}, School of Mathematics and Statistics,
Lanzhou University, Lanzhou 730000, Gansu, P.R. China.

\noindent
E-mail address: \texttt{lihg2024@lzu.edu.cn}
\end{document}